\documentclass[11pt]{amsart}

\usepackage {amscd}

\newtheorem{lemma}{Lemma}

\newtheorem*{conj}{Conjecture}
\newtheorem*{hoved}{Main Theorem}

\newcommand{\pN}{{{\bf P}^N}}
\newcommand{\pk}{{{\bf P}^k}}
\newcommand{\op}{{\mathcal O}}
\newcommand{\gA}{{\mathcal A}}
\newcommand{\gB}{{\mathcal B}}
\newcommand{\gC}{{\mathcal C}}

\newcommand{\ptre}{{{\bf P}^3}}
\newcommand{\pfir}{{{\bf P}^4}}
\newcommand{\pr}[1]{{{\bf P}^{#1}}}
\newcommand{\opN}{\op_\pN}
\newcommand{\opk}{\op_\pk}
\newcommand{\opt}{\op_\ptre}
\newcommand{\opf}{\op_\pfir}

\newcommand{\cod}{\mbox{cod}\,}
\newcommand{\cok}{\mbox{cok}\,}
\newcommand{\rk}{\mbox{rk}\,}
\newcommand{\bt}{{\bf t}}
\newcommand{\bv}{{\bf v}}

\newcommand{\bm}{{\bf m}}

\newcommand{\lpd}{\mbox{lpd}\,}
\newcommand{\tp}{\tilde{p}}
\newcommand{\tq}{\tilde{q}}
\newcommand{\hele}{{\bf Z}}

\newcommand{\te}{\otimes}
\newcommand{\sus}{\subseteq}
\newcommand{\lpil}{\longrightarrow}
\newcommand{\pil}{\rightarrow}
\newcommand{\avb}[1]{\stackrel{#1}\longrightarrow}

\newcommand{\gI}{{\mathcal I}}
\newcommand{\gF}{{\mathcal F}}
\newcommand{\gE}{{\mathcal E}}
\newcommand{\gO}{{\mathcal O}}

\newcommand{\dt}{{\displaystyle \cdot}}

\begin{document}

\title {Monads on projective space}
\author { Gunnar Fl{\o}ystad}
\address{ Matematisk Institutt\\
          Johs. Brunsgt. 12 \\
          5008 Bergen \\
          Norway}   
        
\email{ gunnar@mi.uib.no }

\maketitle

\section*{Introduction}

A monad on projective $k$-space $\pr{k}$ over a field $K$ is a complex
\[ \gA \avb{\alpha} \gB \avb{\beta} \gC \]
of vector bundles on $\pr{k}$ where $\alpha$ is injective and
$\beta$ is surjective.
In this paper we classify completely when there exists monads on
$\pr{k}$ whose maps are matrices of linear forms, i.e. monads
\[  \opk (-1)^a \avb{\alpha}  \opk{}^b \avb{\beta}  \opk (1)^c.  \]
We prove the following.

\begin{hoved} Let $k \geq 1$. There exists a monad as above if and only
if at least one of the following holds.
\begin{enumerate}
\item $ b \geq 2c+k-1$ and $b \geq a+c$. 
\item $ b \geq a+c+k$. 
\end{enumerate}

If so, there actually exists a monad with the map 
$\alpha$ degenerating in expected codimension $b-a-c+1$.
\end{hoved}

Thus in case 2 one sees that there exists a monad whose cohomology is a sheaf 
of constant rank $b-a-c$, i.e. a vector bundle.

\medskip

Our way to this result came through work on curves in $\ptre$. Here monads of
the form
\[ \opt(-1)^{n+r-2} \avb{\alpha} \opt{}^{2n+r-1} \avb{\beta} \opt(1)^n \]
naturally occur. The monad has cohomology $\gI_C(r-2)$, 
the ideal sheaf of a space
curve $C$ twisted with $r-2$. It was known that these monads existed for all
$r \geq 3$ (with $C$ a smooth curve in fact). But it is fairly easy to see
geometrically that the monad could not exist on $\ptre$ for $r \leq 2$.

Consider now the monad
\begin{equation} \opk(-1)^{n+r-2} \avb{\alpha} \opk{}^{2n+r-1} \avb{\beta} \opk(1)^n. 
\label{krn}
\end{equation}
When $k=2$ it is easily seen that it exists for $r \geq 2$ but not for 
$r \leq 1$. When $k=4$ considerations also indicated that it existed for 
$r \geq 4$ but not for $r \leq 3$. Whence we were lead to formulate the 
above result.

\medskip

An interesting geometric consequence of the above result may be obtained by 
letting $k=r=4$ in (\ref{krn}). 
By the theorem above, $\alpha$ may be assumed to degenerate
in codimension 2 along a (locally Cohen-Macaulay) surface $S \sus
\pr{4}$. 
The cohomology of the monad 
is then $\gI_S(2)$. By standard exact sequences one easily sees that
$H^1 \gI_S(1) = n$. Thus $H^0 \gO_S (1) = n+5$ for $n \geq 0$.
This means that $S$ embeds linearly normally into $\pr{n+4}$ and that it 
projects isomorphically down to $\pr{4}$. By a classical theorem of Severi
it is known that the only {\it smooth} 
surface in $\pr{n+4}$ where $n \geq 1$ that enjoys this property is
the Veronese surface
in $\pr{5}$. In case $n=1$ the surface $S$ above is a degeneration of 
the Veronese surface.


\medskip

We would like to motivate the theorem in a more general context.
Let $\gF$ be a sheaf on projective space $\pk$. Suppose $k-r = \lpd \gF$, the
local projective dimension of $\gF$. (If for instance $\gF = \gI_X$,
the ideal sheaf of a smooth projective variety 
$X \sus \pk$ with $\dim X = r$, this holds.)
Then there are, [Wa] Proposition 1.3,  {\it canonical} 
complexes

\[ \gF_i^{\dt} \,: \,\, 0 \pil \gF_i^{-(k-1-i)} \pil \gF_i^{-(k-2-i)} \pil \cdots
\pil \gF_i^0 \pil \cdots \pil \gF_i^i \pil 0 \]
for $i = 0, \ldots, r$ with the following properties :
Each $\gF_i^j$ is a finite sum of line bundles. The cohomology 
$H^j(\gF_i^{\dt}) = 0$ for $j \neq 0$ and $H^0 (\gF_i^{\dt}) = \gF$.
(When $i=0$ this is the sheafification of a minimal resolution
of the graded $K[x_0, \ldots, x_k]$-module 
$\oplus_{n \in \hele} \Gamma(\pk, \gF(n))$.)
Classification of when such complexes exist might be a 
way to understand what kind of algebraic or geometric objects which can exist
on a projective space. The theorem above may be seen as a very small 
contribution to this. 

\medskip

For a more specific motivation consider the complex

\[ \gF_2^{\dt}\, : \,\, 0 \pil \opf(-1)^{n+r-1} \pil \opf{}^{2n+m+r+1}
                   \pil \opf(1)^{n+2m} \pil \opf(2)^m \pil 0. \]
We let the grading of the complex correspond to the twist of the line bundles.
Suppose the complex has cohomology only in degree 0. The cohomology will then
be a sheaf $\gE$ of rank $2$. If this is to be a vector bundle one must have
$c_3(\gE) = 0$ and $c_4(\gE) = 0$. This gives two equations relating
$m,n$ and $r$. 
The following integer values of $r$ (and maybe others) will give integer
solutions for $m$ and $n$.
\begin{itemize}
\item[a.] For $r=1$ one gets $(m,n) = (0,0)$ or $(m,n) = (0, -1)$.
When $(m,n) = (0,0)$ one gets the vector bundle
$\opf{}^2$. When $(m,n) = (0,-1)$ one gets the bundle $\opf(1) \oplus \opf(-1)$.
Here one has to shift a term across an arrow if the summand is of negative 
order.
\item[b.] 
For $r=0$ one gets $(m,n) = (0,0)$ or $(m,n) = (2,6)$.
When $(m,n) = (0,0)$ one gets the bundle $\opf \oplus \opf(-1)$. 
When $(m,n) = (2,6)$ one actually gets the $\gF_2^{\dt}$ complex of the
Horrocks-Mumford bundle. 

\item[c.] For $r= -1$ one gets $(m,n) = (0,0)$ or $(m,n) = (11,22)$.
When $(m,n) = (0,0)$ one gets the bundle $\opf(-1)^2$. When
$(m,n) = (11,22)$ one might get a potentially {\it new rank 2 bundle} 
$\gE$ on $\pr{4}$ with
$\gF_2^{\dt}$ complex
\[ 0 \lpil \opf(-1)^{20} \lpil \opf{}^{55} \lpil \opf(1)^{44} 
\lpil \opf(2)^{11} \lpil 0. \]
The bundle $\gE$ would have $c_1(\gE) = -2$ and $c_2(\gE) = 12$.
\end{itemize}
\medskip 

The next step in classifying complexes on projective space might be to
classify complexes
\[ 0 \lpil \opk(-2)^c \avb{\beta^{\vee}} \opk(-1)^b \avb{\alpha^{\vee}}
           \opk{}^a  \]
where the cokernel in degree $0$ is $\gI_X(2c-b)$ for a locally Cohen-Macaulay
subscheme $X \sus \pk$ of codimension $2$. This means that $\beta$ is 
surjective and $\alpha$ degenerates in codimension $2$. Note that we must have 
the relation
$b+1 = a+c$. (The Main Theorem deals with the case $a+c \leq b$.) 
We may then write the complex as
\[ \opk(-2)^n \avb{\beta^{\vee}} \opk(-1)^{2n+r} \avb{\alpha^{\vee}} 
   \opk{}^{n+r+1}. \]

\begin{conj} Let $k \geq 2$. 
The complex above with $\beta$ surjective and $\alpha$
degenerating in codimension $2$ exists if and only if $r \geq 0$ and
$ n \leq \binom{r+3-k}{2}$.
\end{conj}

For instance if $k=2$ this is readily seen to hold. When $n = \binom{r+1}{2}$
and $S = K[x_0,x_1,x_2]$ the complex above corresponds to the resolution of 
the power $\bm^r \sus S$ where $\bm = (x_0, x_1, x_2)$.

When $k=3$ the necessity of the above condition on $n$ is also readily seen 
geometrically to hold.
\vskip 3mm
Let us just say some words about the organization of the paper. In the
first section we give the existence of the monads. This is based on
an immediate explicit construction. The main work of the section is to 
verify the statement about the degeneracy locus.

In the second section we show the necessity of the conditions given.
This is an argument based on the study of degeneracy loci.

We work over an arbitrary field $K$.

\section{ Existence }

We shall prove the existence of the monads by providing an
explicit construction. Denote by $X_{r,r+n}$ the $r$ by $r+n$ matrix
\[ \left(  \begin{array}{ccc ccc c} 
          x_0 & x_1 & \ldots & x_n & & &  \\
              & x_0 & x_1 & \ldots & x_n &  &  \\
              &     & \ddots & & & \ddots & \\
              &     &        & x_0 & x_1 & \ldots & x_n  
\end{array} \right). \]

A basic fact we may note is that $X_{r, r+n}$ degenerates in rank if and only
if all $x_i = 0$. Similarly denote by $Y_{r, r+m}$ the $r$ by $r+m$ matrix
\[ \left(   \begin{array}{ccc ccc c} 
          y_0 & y_1 & \ldots & y_m & & &  \\
              & y_0 & y_1 & \ldots & y_m &  &  \\
              &     & \ddots & & & \ddots & \\
              &     &        & y_0 & y_1 & \ldots & y_m  
\end{array}  \right).  \]

Let $\sigma_k = \sum_{k = i + j} x_i y_j$. Form the $r$ by $r+m+n$
matrix 

\[ \Sigma_{r,r+m+n} = 
\left(   \begin{array}{ccc ccc c} 
     \sigma_0 & \sigma_1 & \ldots & \sigma_{m+n} & & &  \\
              & \sigma_0 & \sigma_1 & \ldots & \sigma_{m+n} &  &  \\
              &     & \ddots & & & \ddots & \\
              &     &        & \sigma_0 & \sigma_1 & \ldots & \sigma_{m+n}  
\end{array} \right).  \]

The following lemma is easily verified.

\begin{lemma}
\[ X_{r,r+n}\cdot Y_{r+n,r+n+m} = Y_{r, r+m} \cdot X_{r+m,r+m+n} = 
\Sigma_{r,r+n+m}. \qed \]
\end{lemma}

\noindent {\sl Note.} We have been notified that V. Ancona and G. Ottaviani
 used the same matrices and
lemma (with $n=m$) in [An-Ot].

\medskip

Let 
\[ S = K[x_0, \ldots, x_n, y_0, \ldots y_m] \]
and $\pN = \mbox{Proj}\, S$ where $N = n+m+1$.
We may form the complex
\begin{equation}
 \op_{\pN}(-1)^{r+n+m} \avb{\alpha} \op_{\pN}{}^{2r+n+m} \avb{\beta} 
     \op_{\pN}(1)^r, \label{clx}
\end{equation}
where the maps $\beta$ and $\alpha$ are given by the matrices
\[ B = \left [ \begin{array} {cc}
               X_{r,r+n}  & Y_{r,r+m} 
               \end{array} \right ] , \quad
   A = \left [ \begin{array} {c}
               Y_{r+n,r+n+m} \\
               -X_{r+m, r+m+n} 
                \end{array} \right ]. \]

We easily see that the map $\beta$ is surjective. The whole existence part
of the Main Theorem may be derived from this complex 
as we shall shortly see, but 
first we investigate the degeneracy loci of the map $\alpha$.

\medskip

Let $Z_d \sus \pN$ be the locus where $A$ degenerates to rank $r+n+m-d$.

\begin{lemma}

\begin{enumerate}
\item If $d > \max(n,m)$ then $Z_d = \emptyset$.
\item If $\min(n,m) < d \leq \max(n,m)$ then 
$d = \max(n,m)$ and $\cod Z_d = \min(n,m) + 1$.
\item If $d \leq \min(n,m) + 1$ then $\cod Z_d \geq d$.
\end{enumerate}
\end{lemma}

In particular we see that if $|n-m| \leq 1$ then $\cod Z_d \geq d$ for
all $d \geq 0$.
\medskip
\begin{proof} 1. Suppose $d > \max(n,m)$. If $\rk A \leq r+n+m-d$ then
$\rk A < r+n$ and $\rk A < r+m$. Then it is easily seen that all $x_i =0$ and
all $y_i = 0$.

2. In this case assume $m < n$. 
Since $d > m$ we get $\rk A \leq r+n+m-d < r+n$. This gives that all
$y_i = 0$. Since $d \leq n$ we also see that $\rk A = r+n+m-d \geq r+m$.
We see that we must have $\rk A = r+m$ and $d=n$.


3. We will show that there is a linear subspace of $\pN$ of dimension
$d-1$ where $A$ does not degenerate to rank $r+n+m-d$. This will prove
the third part of the lemma.

Let $x_i = t_i$ for $i = 0,\ldots, d-1$ and $x_i = 0$ for $i \geq d$.
Let $y_{m-i} = t_{d-1-i}$ for $i = 0,\ldots, d-1$ and $y_{m-i} = 0$
for $i \geq d$.
This gives a linear subspace $L_{d-1} \sus \pN$ of dimension $d-1$.
With these substitutions the matrix $A$ now takes a form (letting
$\bt = (t_0, \ldots, t_{d-1}))$
\[ A_\bt = \left [ \begin{array} {c}
                   T_1 \\     
                   -T_2  \end{array} \right ], \]
where 
\[ T_1 = \left(   \begin{array}{ccc ccc ccc} 
          0 & \ldots & 0 & t_0 & \ldots &t_{d-1} &&&  \\
              & 0 & \ldots & 0 & t_0 & \ldots & t_{d-1} & &\\
              &     & \ddots & & &\ddots &  & \ddots& \\
              &     &        & 0 & \ldots & 0 & t_0 & \ldots & t_{d-1}   
\end{array}  \right)  \] 
and
\[ T_2 = \left(   \begin{array}{ccc ccc ccc} 
         t_{0} & \ldots & t_{d-1} & 0 & \ldots &0 &&&  \\
              & t_{0} & \ldots & t_{d-1} & 0 & \ldots & 0 & &\\
              &     & \ddots & &\ddots & & & \ddots& \\
              &     &        & t_{0} & \ldots & t_{d-1} & 0 & \ldots & 0   
\end{array}  \right)  \]

If $d=1$ then clearly $\rk A = r+n+m$, so suppose $d > 1$.
If $t_{d-1} = 0$ then a submatrix of $A_\bt$ is of the form
\[ A^\prime_{\bt^\prime} = \left [ \begin{array} {c}
                   T_1^\prime \\     
                   -T_2^\prime  \end{array} \right ], \]
which is the matrix corresponding to the case $n^\prime = n$,
$m^\prime = m-1$ and $\bt^\prime = (t_0, \ldots, t_{d-2})$. 
By induction $\rk A^\prime_{\bt^\prime} > r+n^\prime + m^\prime - (d-1) =
r+n+m-d$.

In the following for a vector $\bv = (v_1, \ldots, v_{r+n+m})$ in
$K^{r+n+m}$, let
\[  \min(\bv) = \min\{ i \, | \, v_i \neq 0\}. \]
Assume now that there exists
a $\bt^0$ with $t_{d-1}^0\neq 0$ such that $\rk A_{\bt^0} \leq r+n+m-d$.
Then there is a $d$-dimensional subspace
of $\bv$ in $K^{r+n+m}$ such that $A_{\bt^0} \cdot \bv^t = 0$. Among these
there must be a non-zero $\bv$ such that $\min(\bv) \geq d$. But from
the equation 
$ T_2 \cdot \bv^t = 0 $
and the fact that $t^0_{d-1} \neq 0$ we see that we must in fact have
$\min(\bv) \geq r+m+d$. In particular $\min(\bv) \geq m + 1$. 
But this is impossible since we have $T_1 \cdot \bv^t = 0$ 
with $t_{d-1}^0\neq 0$.
\end{proof}

Let 
\[ \opN(-1)^{r+n+m-s} \avb{\phi} \opN(-1)^{r+n+m}\] 
be a general injection. Recall the map $\alpha$ from (\ref{clx}).

\begin{lemma}
Suppose $\cod Z_d \geq d$ for $d = 1, \ldots, s+1$ and $\cod Z_d \geq
s+1$ for $ d \geq s+1$. Then $\alpha \circ \phi $
degenerates in codimension $s+1$.
\end{lemma}

In particular we see that if $|n-m| \leq 1$ then $\alpha \circ \phi$ will 
degenerate in codimension $s+1$.

\begin{proof} Let $\gE = \ker \beta$.
Dualizing 
\[ \op_{\pN}(-1)^{r+n+m} \avb{\alpha} \gE \]
we get
\[ \gE^{\vee}(-1) \avb {\alpha^{\vee}(-1)} \op_{\pN}{}^{r+n+m}. \]
Let 
\[ T_0 = \cok \alpha^{\vee}(-1). \]
Then $T_0$ is generated by the sections coming from $\gO_{\pN}{}^{r+n+m}$.
We get for $T_0$ a flattening stratification of locally closed
subschemes (see \cite{Mu} section 8 or \cite{Ei} section 20.2)
\[ \pN = \cup_{d \geq 0} Z_{d,0} \]
such that $T_{0|Z_{d,0}}$ is locally free of rank $d$.
Note that $Z_{d,0} = Z_d$. 

Take a general section of $\op_\pN{}^{r+n+m}$. 
This gives a section of $T_0$ and
a sequence
\[ \op_\pN \lpil T_0 \lpil T_1 \lpil 0. \]
Let 
\[ \pN = \cup_{d \geq 0} Z_{d,1} \] be a flattening stratification for $T_1$.
We must have 
\[ Z_{d,1} \sus Z_{d,0} \cup Z_{d+1,0}. \]
Since the sections of $\op_\pN{}^{r+n+m}$ generate $T_{0|Z_{d,0}}$, 
we see that
\[ \cod (Z_{d,1}\cap Z_{d,0}) > \cod Z_{d,0} \]
for $d \geq 1$. Thus in this case
\[ \cod Z_{d,1} \geq \cod Z_{d,0} + 1 \]
or 
\[ \cod Z_{d,1} \geq \cod Z_{d+1,0}. \]

In this way we may proceed to $T_s$. This sheaf has again a flattening
stratification 
\[ \pN = \cup_{d \geq 0} Z_{d,s} \] 
such that when $d \geq 1$ we have
\[ \cod Z_{d,s} \geq \cod Z_{d+a,0} + s - a \]
for some $a$ with  $ 0  \leq a  \leq s$. 
By hypothesis we get 
$\cod Z_{d,s} \geq s+1$ for $d \geq 1$.

The above process gives a diagram
\[ \begin{CD} 
   @.  \opN^s  @=  \opN{}^s \\
   @.   @VVV  @VVV \\
  \gE^{\vee}(-1)  @>>> \opN{}^{r+n+m}  @>>>  T_0 \\
   @|   @VVV  @VVV \\     
      \gE^{\vee}(-1) @>>>   \opN{}^{r+n+m-s}  @>>>  T_s
\end{CD} \]

We may take $\phi^{\vee}(-1)$ to be the lower middle vertical map.
Then $\alpha \circ \phi$ degenerates in 
$\mbox{Supp}\, T_s \sus \cup_{d \geq 1} Z_{d,s}$
which has codimension $\geq s+1$. But since $s$ is the difference of
ranks between $\gE^{\vee}(-1)$ and $\opN{}^{r+n+m-s}$, we know that the 
codimension is at most $s+1$. Thus it is exactly $s+1$.
\end{proof}

\begin{proof} [Proof of the existence part of the Main Theorem.]
Given $a,b,c$ and $k$ satisfying part 1 of the Main Theorem.
Let $r=c$. Choose $n$ and $m$ with $|n-m| \leq 1$ such that $b-2c = n+m$.
By (\ref{clx}) there exists a monad
\[ \opN (-1)^a \avb{\alpha} \opN{}^b \avb{\beta} \opN{}^c \]
with $N = n + m + 1$. By Lemmata 2 and 3 we may assume that $\alpha$ 
degenerates in codimension $b-c-a+1$. 
By restricting to a general subspace $\pk \sus \pN $ we get part 1 of the
Main Theorem.

To prove the existence in part 2 we may assume that $b \leq 2c+k-1$ since 
else we may refer to part 1. 
Since $b \geq a+c+k $ we get $a \leq c-1$. Thus
\[ b \geq a+c+k \geq 2a+1+k \geq 2a+k-1. \]
But then by part 1 there exists a monad
\[ \opk(-1)^c \avb{\beta^{\vee}} \opk{}^b \avb{\alpha^{\vee}} \opk(1)^a \]
with $\beta^{\vee}$ degenerating in codimension $b-a-c+1 \geq k+1$.
But then $\beta^{\vee}$ does not degenerate. Dualizing we get a monad
\[  \opk(-1)^a \avb{\alpha}  \opk{}^b \avb{\beta}  \opk(1)^c \]
with $\alpha$ not degenerating. 
\end{proof}

\section{ Necessity of conditions}

Suppose now we have given a monad
\[ \opk(-1)^a \avb{\alpha} \opk{}^b \avb{\beta} \opk(1)^c . \]
We wish to prove the numerical conditions on $a, b, c$ and $k$ given in the
Main Theorem. The image of 
\[ \Gamma( \opk{}^b) \avb{\Gamma(\beta)} \Gamma (\opk(1)^c) \]
determines a subspace $V \sus \Gamma( \opk(1)^c) $
which generates the bundle $\opk(1)^c$ since $\beta$ is surjective.
Also $\dim V \geq c+k$ since otherwise $\beta$ would degenerate in a non-empty 
subscheme of codimension $\dim V - c + 1$ by \cite{Fu} 14.4.13.
Let $U \sus V $ be a general subspace of dimension $c+ k -1$.
Then the map
\[ U \te \opk \lpil \opk(1)^c \]
degenerates in dimension $0$, again by \cite{Fu} 14.4.13.
Fix a splitting
\[ \Gamma(\opk{}^b) \avb{\longleftarrow} V. \]
Let $W = \Gamma(\opk{}^b)/ U $ and $S = K[x_0, \ldots, x_k]$. We get a diagram of
free $S$-modules.
\[ 
\begin{CD}
     @. U \te S @= U \te S  \\
       @.              @VVV     @VVpV \\
    S(-1)^a @>>> S^b @>>> S(1)^c\\
     @| @VVV  @. \\
   S(-1)^a @>q>>  W \te S  @.
\end{CD} \]
Let $\tp$ and $\tq$ denote the corresponding maps of sheaves. We note that
there is a surjection 
\[ \cok \tq \lpil \cok \tp \lpil 0. \]
Since $\tp$ degenerates in expected codimension, by [Bu-Ei] Theorem 2.3
we have an equality
\[ \mbox{Fitt}_1(\cok \tp) = \mbox{Ann} (\cok \tp) \]
where $\mbox{Fitt}_1(\cok \tp)$ is the first Fitting ideal generated by the 
$c \times c$ minors of the matrices (locally) representing $\tp$. We now get
\[ \mbox{Fitt}_1(\cok \tq) \sus \mbox{Ann}( \cok \tq) 
   \sus \mbox{Ann}(\cok \tp) = \mbox{Fitt}_1(\cok \tp) \]
where the first inclusion is valid if we replace $\cok \tq$ by any 
coherent sheaf, \cite{Ei} Proposition 20.7.a.
This gives
\[ \mbox{Fitt}_1(\cok q) \sus \Gamma_* \mbox{Fitt}_1(\cok \tq) 
   \sus \Gamma_* \mbox{Fitt}_1(\cok \tp). \]
Since $p$ degenerates in expected codimension $k$, and $S$ is a Cohen-Macaulay
ring, $ S / \mbox{Fitt}_1( \cok p)$
will be a Cohen-Macaulay ring of dimension 1 by \cite{Ei} Theorem 18.18 or
\cite{Fu} Theorem 14.4.c.
Thus the irrelevant maximal ideal $\bm \sus S$ is not an associated
prime of $\mbox{Fitt}_1(\cok p)$ and this is thus a saturated ideal. This gives
\[ \Gamma_* \mbox{Fitt}_1(\cok \tp) = \mbox{Fitt}_1(\cok p). \]
Since now $\mbox{Fitt}_1(\cok p)$ is generated by polynomials of degree
$\geq c$, no polynomial in $\mbox{Fitt}_1(\cok q)$ will have degree $< c$.
Note that since $\alpha$ is injective and $S^b \pil W \te S$ 
is a general quotient, the map 
$q$ may be assumed to generically have maximal rank.
If therefore $q$ is generically surjective, we must have 
\[ \dim W \geq c. \]
Otherwise we have  
\[ \dim W >  a. \]
Since $\dim W = b-c-k+1$ this gives
\[ b \geq 2c+k-1 \]
or \[ b \geq a+c+k. \] 
This proves the necessity of the conditions in the Main Theorem.

\end{document}